\documentclass[review]{elsarticle}

\usepackage{lineno,hyperref}
\usepackage{graphicx}
\usepackage{float}
\usepackage{xcolor}

\modulolinenumbers[5]

\journal{arXiv}

\usepackage{amsmath,amssymb}
\usepackage{amsthm}
\usepackage{todonotes,commath}

\usepackage{algorithmic}

\newtheorem{thm}{Theorem}
\newtheorem{cor}[thm]{Corollary}
\newtheorem{lem}[thm]{Lemma}
\renewcommand{\L}{\mathcal{L}}
\newcommand{\e}{\mathrm{e}}

\DeclareMathOperator{\PR}{\mathbb{P}}

\DeclareMathOperator{\Var}{\mathrm{Var}}

\begin{document}
\author[Durham]{Robert~E.~Gallagher}
\author[Durham,Turing]{Louis~J.~M.~Aslett}
\author[Oxford]{David Steinsaltz}
\author[WashU,Turing]{Ryan~R.~Christ}

\begin{frontmatter}

\title{Improved Concentration Bounds for Gaussian Quadratic Forms}

\address[Durham]{Department of Mathematical Sciences, Durham University, United Kingdom}
\address[Turing]{The Alan Turing Institute, London, United Kingdom}
\address[Oxford]{Department of Statistics, Oxford University, United Kingdom}
\address[WashU]{McDonnell Genome Institute, Washington University in St. Louis, United States}

\begin{abstract}
 For a wide class of monotonic functions $f$, we develop a Chernoff-style concentration inequality for quadratic forms $Q_f \sim \sum\limits_{i=1}^n f(\eta_i) (Z_i + \delta_i)^2$, where $Z_i \sim N(0,1)$.  The inequality is expressed in terms of traces that are rapid to compute, making it useful for bounding p-values in high-dimensional screening applications.  The bounds we obtain are significantly tighter than those that have been previously developed, which we illustrate with numerical examples.  
\end{abstract}

\begin{keyword}
quadratic form \sep
generalized non-central chi-square distribution \sep
concentration inequality \sep
Hilbert-Schmidt Information Criteria \sep
tail bound
\end{keyword}

\end{frontmatter}

\linenumbers

\section{Introduction and Background}

We consider the problem of finding an upper bound for the cumulative distribution function (cdf) 
of random variables of the form $Q_f \sim \sum\limits_{i=1}^n f(\eta_i) (Z_i + \delta_i)^2$, where $Z_i \sim N(0,1)$, $f: \mathbb{R} \rightarrow \mathbb{R}$, and $\delta_i$ and $\eta_i$ are deterministic scalars. Many applications lead to this form with $\{\eta_i\}_{i=1}^n$ being the eigenvalues of a symmetric matrix $M \in \mathbb{R}^{n \times n}$; for example, a quadratic form $X^\top f(M) X$ where $X \sim N(\mu,I)$ and $f(M)$ represents $f$ applied to the eigenvalues of $M$.  As described in \citet{Christ:2017aa} and \citet{Christ2020}, results of this kind can be generalized to cases where $M$ is asymmetric with careful treatment of $f$.  

$Q_f$ arises as the limiting distribution of test statistics used in a wide range of applications.  These statistics include the Hilbert-Schmidt Information Criterion used for high-dimensional independence testing \cite{gretton2005, Zhang2018}, score statistics for linear and genearlized linear mixed models commonly used in genomics \cite{lin1997variance, wu2011rare}, and the goodness-of-fit statistic proposed by \citet{pena2002powerful} for ARMA models in time series analysis. It is easy to see that $Q_f$ has mean
\[
\mathbb{E}(Q_f) = \sum\limits_{i=1}^n  f(\eta_i) \delta_i^2 + \sum\limits_{i=1}^n  f(\eta_i)
\]
and variance
\[
\Var(Q_f) = 2 \left( \sum\limits_{i=1}^n  f(\eta_i)^2 \delta_i^2 + 2 \sum\limits_{i=1}^n f(\eta_i)^2\right).\]
Work in \cite{Christ2020} established a concentration inequality to bound the tails of $Q_f$, which yield a set of bounds for different functions. The results of \cite{Christ2020} show that it is possible to find polynomial bounds,
but these are not constructed explicitly.  We provide here explicit optimal coefficients for bounds of this form in the single-spectrum case. 
This earlier work yielded the following bound on $Q$ (by which we designate the base version of
$Q_f$, where $f$ is the identity function):

\begin{thm}[see p.75 in \cite{Christ:2017aa}]
\label{theorem1}
Let $X \sim N(\mu,I)$ and $M$ be a real, symmetric matrix. Let $Q=X^\top M X$. Let $\nu = 2 \left( 4 \left| \left| M \mu \right| \right|^2_{2} + 2 \left| \left| M \right| \right|^2_{HS} \right)$ and let $b=\underset{i}{\max}\left| \lambda_i\right|$, where $\{\lambda_i\}_{i=1}^n$ are the eigenvalues of $M$.  Then, for all $q > \mathbb{E}\left[ Q \right]$, 

\begin{equation*}
\mathbb{P}\left( Q  > q \right) \leq  \left\{ \begin{array}{lcc}
 \exp\left( -\frac{1}{2} \frac{(q-\mathbb{E}\left[ Q \right])^2}{\nu}\right) &  & \mathbb{E}\left[ Q \right] < q  \leq \frac{\nu}{4 b}+ \mathbb{E}\left[ Q \right] \\
 \exp\left( \frac{1}{2} \frac{\nu}{\left(4b\right)^2 } \right) \exp\left( -\frac{q-\mathbb{E}\left[ Q \right]}{4b}\right) & & q > \frac{\nu}{4 b }+\mathbb{E}\left[ Q \right].
\end{array}\right.
\end{equation*}
Similarly, for all $q < \mathbb{E}\left[ Q \right]$,
\begin{equation*}
\mathbb{P}\left( Q  < q \right) \leq  \left\{ \begin{array}{lcc}
 \exp\left( -\frac{1}{2} \frac{(\mathbb{E}\left[ Q \right]-q)^2}{\nu}\right) &  & \mathbb{E}\left[ Q \right] -\frac{\nu}{4 b} \leq q  < \mathbb{E}\left[ Q \right] \\
 \exp\left( \frac{1}{2} \frac{\nu}{\left(4b\right)^2 } \right) \exp\left( -\frac{\mathbb{E}\left[ Q \right]-q}{4b}\right) & & q < \mathbb{E}\left[ Q \right]-\frac{\nu}{4 b }.
\end{array}\right.
\end{equation*}
\vspace{1em}
\end{thm}

The proof of this result relies on a Chernoff-style bound involving the cumulant generating function (cgf) of $Q$, which has two main types of terms:
\begin{align*}
	\L_1(x)&=-\log(1-2x)/2 \quad\text{ and}\\
	\L_2(x)&=\frac{x}{1-2x}.
\end{align*}
Each of these is bounded by a quadratic function, leading to an overall bound in terms of easily computable
coefficients.  We improve on this previous work by constructing a family of quadratics that yield pointwise tighter bounds on $\L_1$ and $\L_2$.  We then show how these can be incorporated into an optimisation step to yield tighter bounds on the tails of $Q_f$.

In Section 2 we present our main results.  First we present Lemma \ref{lemma1}, which tightens the quadratic bounds above from \cite{Christ:2017aa}.  From this lemma, we derive the corresponding improved bounds on the tails of $Q_f$ in Theorem \ref{Thm1}. 
Specialisation of these results to some particular functions $f$ then follow in corollaries. 
In Section 3, we empirically demonstrate the improvement provided by these bounds with an application to a simulatd matrix with a exponentially decaying spectrum.  Section 4 concludes with discussion of potential future improvements. Proofs for the main results are presented in Section 5.

\section{Main Results}


Our results depend upon elementary upper bounds on $\L_1(x)$ and $\L_2(x)$ in the form of parabolas passing through the origin. We describe the coefficients of these parabolas 
in terms of the width of the (symmetric) interval on which the bounds are to be applied, and on the parameter $t$
that arises from the cgf. We exploit two openings for improvement: optimising the coefficients of the parabola and optimising the width of the scaled domain over which it bounds $\L_1 \left(tf(x)\right)$ and $\L_2 \left(tf(x)\right)$.
\begin{lem}
\label{lemma1}
Let $f(x)$ be a monotonic increasing function such that $f(0)=0$. Let $L$ be a fixed positive real number, and $t \in [0,t^{\star})$, where
\[
  t^{\star} = \min\bigl\{|1/2f(L)|,|1/2f(-L)| \bigr\}.
\]
Furthermore, suppose that over the region $x \in (0,L], t \in [0, t^{\star})$ 
the following inequalities are satisfied for both $\L_1\left(tf(x)\right)$ and $\L_2 \left(tf(x)\right)$:
\begin{align}
	x\bigl(\partial_x \L\left(tf(x)\right)/2 + t f'(0) \bigr) &\geq \L\left(tf(x)\right), \label{E:lemcond1}\\
	\frac{\L\left(tf(x)\right) - \L\left(tf(x)\right)}{2x} &\geq t f'(0). \label{E:lemcond2}
\end{align}
For each $t \in [0,t^{\star})$ define
\begin{align*}
\alpha_f (L,t) &= \L_1\left(tf(x)\right)/ L^2 - t f'(0)/L, \\
	\beta_f (L,t) &=   \L_2 \left(tf(x)\right) / L^2 - t f'(0)/L,\\
\gamma_f(t) &= t f'(0) .
\end{align*}

Then for each $t \in [0,t^{\star})$,
among all quadratic function $x \mapsto ax^2 + b x$ 
that maintain $g^1_t(x)\le 0$ over the whole region $x \in [-L,L]$, where
\[
	g^1_t(x) := \L_1 \left(tf(x)\right) - \bigl( a x^2 + b x \bigr),
\]
the difference $|g^1_t(x)|$ is minimised at every point $x$
by the choice $a= \alpha_f(L,t)$  and $b=\gamma_f(t)$;
and among those that maintain $g^2_t(x)\le 0$ over the whole region $x \in [-L,L]$, where
\[
g^2_t(x) := \L_2 \left(tf(x)\right) - \bigl( a x^2 + b x \bigr),
\]
the difference $|g^2_t(x)|$ is minimised at every point $x$
by the choice $a= \beta_f(L,t)$ and $b=\gamma_f(t)$.
\end{lem}

This lemma will allow us to build on the existing result from \cite{Christ:2017aa}. In the original 
form of this theorem $t$ was restricted so that $tf(x) < 1/4$, avoiding the asymptote at $1/2$. We remove this boundary at 1/4 and allow $tf(x)$ to get arbitrarily close to $1/2$. We also reinterpret $L$, so that it now defines the 
domain of $x$ rather than that of $tf(x)$. It also means that for every endpoint along the interval $[-L,L]$ we can obtain optimal coefficients on our quadratic bounds. This yields a new bound on the tails of $Q_f$ as follows.

\begin{thm}
\label{Thm1}
Let $\xi = c \left(\sum\limits_{i=1}^n  \eta_i \delta_i^2 + \sum\limits_{i=1}^n  \eta_i\right)$ where $c=f'(0)$, and let $L$ be set to $\underset{i}{\max}\left| \eta_i\right|$.  Suppose $f$ satisfies the conditions in Lemma \ref{lemma1}. 
Then for all $q > \xi$,
\begin{equation} \label{eq:PQf}
\PR(Q_f > q) \leq \min_{t \in (0, 1/2d)} \bigl[ \exp(\nu_f(t)/2 -\left(q -\xi\right) t) \bigr],
\end{equation}
where $d = \underset{i}{\max}\left| f(\eta_i)\right|$, and
\begin{equation}  \label{eq:nuft}
  \nu_f(t)=2 \left( \beta_f(L,t) \sum\limits_{i=1}^n  \eta_i^2 \delta_i^2 + \alpha_f(L,t) \sum\limits_{i=1}^n \eta_i^2\right).
\end{equation}
Furthermore, for all $q < \xi$,
\begin{equation} \label{eq:PQf2}
  \PR(Q_f < q) \leq \min_{t \in (0, 1/2d)} \bigl[\exp(\nu_f(t)/2 -\left(\xi -q\right) t) \bigr].
\end{equation}
\end{thm}

In the central use case, where $Q_f$ arises as $X^\top f(M) X$, we can apply Theorem \ref{Thm1}, 
where $\xi = c \, (\mu^\top M \mu + \mathrm{tr}\left(M\right))$ and \[\nu_f(t) = 2 \left( \beta_f(L,t) \left| \left| M \mu \right| \right|^2_{2} + \alpha_f(L,t) \left| \left| M \right| \right|^2_{HS}\right).\] This allows us to quickly compute 
tight tail bounds on $X^\top f(M) X$.  In the following corollaries we address special cases of $f$. 

\begin{cor}
	\label{identity_corollary}
Let $f(x)=x$. Then the cdf of $Q_f$ is bounded as in equations \eqref{eq:PQf} and \eqref{eq:PQf2} where in equation \eqref{eq:nuft} we set ${\alpha_f(L,t) = \L_1\left(tL\right)/L^2 - t/L}$ and ${\beta_f(L,t) = tL/(L^2(1-2 tL)) - t/L}$.
\end{cor}
\noindent
\textit{Proof:}
Since $|f(L)| = |f(-L)|$, the $t^{\star}$ from Lemma 1 is equal to $1/2f(L)=1/2L$.
The conditions \eqref{E:lemcond1} and \eqref{E:lemcond2} may be written in terms of
the variable $z = tx$, and these inequalities then need to hold for $z \in [0, 1/2)$.
The two conditions for $\L_1$ become
\begin{align*}
\frac{z}{(1-2z)} + 2z &\geq -\log(1-2z) \quad \text{and} \\
-\log(1-2z) + \log(1+2z) &\geq 4z ,
\end{align*}
while the two conditions for $\L_2$ become
\begin{align*}
\frac{z}{(1-2z)^2} + 2z &\geq \frac{2z}{1-2z} \quad \text{and} \\
\frac{z}{1-2z} + \frac{z}{1+2z} &\geq 2z.
\end{align*}
All of these inequalities hold for $z \in (0,1/2)$, and so Lemma 1 holds where $f$ is
the identity function. The result follows by application of Theorem 1. \qed

\begin{cor} \label{C:xp}
Let $f(x)=x^p$ for some positive integer $p \ge 2$. Then the cdf of $Q_f$ is bounded as in equations \eqref{eq:PQf} and \eqref{eq:PQf2} where in equation \eqref{eq:nuft}, \[{\alpha_f(L,t) = \L_1\left(tL^p\right)/L^2} \quad\text{and}\quad \beta_f(L,t) = tL^{p-2}/(1-2 tL^p).\]

\end{cor}
\noindent
\textit{Proof:}
Since $|f(L)|=|f(-L)|$, the $t^{\star}$ from Lemma 1 is equal to $1/2L^p$. We introduce the variable $z = tx^p$ and note that our original region, $x \in [0,L]$ and $t \in [0, 1/2L^p)$, corresponds
to $z \in [0, 1/2)$.\\
Substituting the definitions of $z, \L_1, \L_2$ into condition
\eqref{E:lemcond1} yields
\begin{align*}
\frac{pz}{1-2z} &\geq -\log(1-2z), \\
\frac{pz}{(1-2z)^2} &\geq \frac{2z}{1-2z}.
\end{align*}
The condition \eqref{E:lemcond2} is trivial for even $p$, while for
odd $p$ it becomes
\begin{align*}
	-\log(1-2z) + \log(1+2z) &\geq 0, \\
	\frac{z}{1-2z} + \frac{z}{1+2z} &\geq 2z.
\end{align*}
All of these inequalities hold for $z \in [0,1/2)$ and $p\ge 2$,
so Lemma \ref{lemma1} holds for $f(x) = x^p$. The result follows by application of Theorem \ref{Thm1}. \qed

With essentially the same proof used for Corollary \ref{C:xp}, we can formulate the result of
Theorem \ref{Thm1} for matrix powers.  Note that in following case, $\xi = 0$.
\begin{cor} \label{C:xpMat}
For any positive integer $p \ge 2$, for each $q>0$
\begin{equation}
\PR(X^\top M^p X > q) \leq \min_{t \in (0, 1/2d)}\e^{-qt + \nu_f(t)/2},
\end{equation}
and for  $q < 0$
\begin{equation}
\PR(X^\top M^p X < q) \leq \min_{t \in (0, 1/2d)}\e^{qt + \nu_f(t)/2},
\end{equation}
where $\nu_f(t)$ is defined in \eqref{eq:nuft} and $\alpha_f(L,t) = -\log(1-2t L^p)/2L^2$, $\beta_f(L,t) = tL^{p-2}/(1-2 tL^p)$.
\end{cor}

\section{Examples}

\begin{figure}[H]
	  \includegraphics[width=\linewidth]{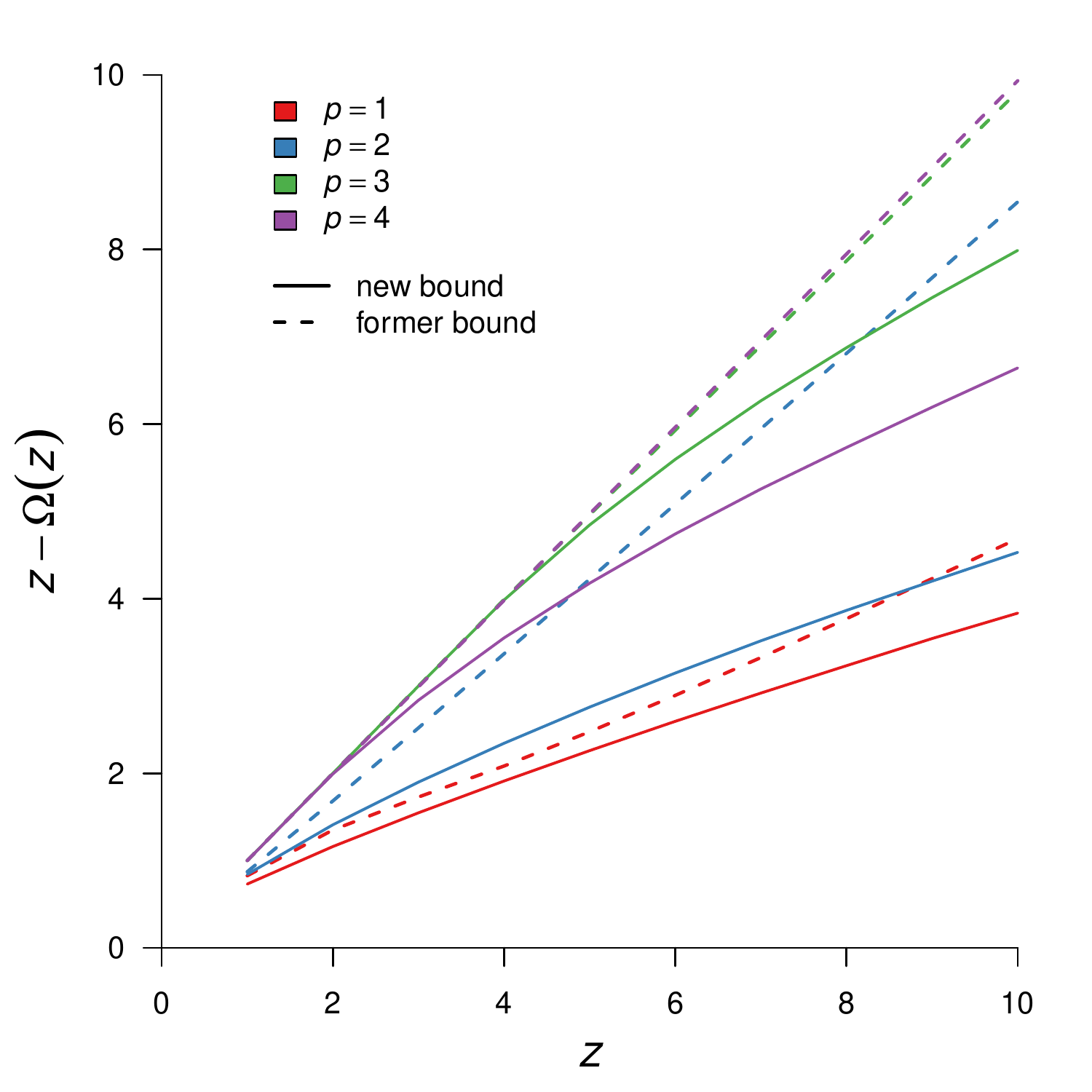}
	\caption{Modified $Q-Q$ plot showing the difference between the negative base 10 logarithm of the true right tail probabilities estimated by our simulations and those estimated by $U_p(q)$.  In other words, we plot $z - \Omega(z)$, where $\Omega\left(z\right) = -\log_{10}\left(1-U_p\left(\hat{F}^{-1}_p\left(1-10^{-z}\right)\right)\right).$
}
	\label{fig:modqq}
\end{figure}

Here we compare the bounds in Corollary \ref{identity_corollary} and Corollary \ref{C:xpMat} to the bounds provided in \citet{Christ:2017aa} and \citet{Christ2020} for different matrix powers $p =1,2,3,4$.  For this comparison, we simluated a matrix with an exponentially decaying spectrum of eigenvalues, a case which is relatively common in applications.  See Figure \ref{fig:modqq}.



\numberwithin{figure}{section}
\noindent For this comparison, we simluated a matrix with an exponentially decaying spectrum of eigenvalues, a case which is relatively common in applications.  See Figure \ref{fig:modqq}.  Note that we have plotted the logarithm (base 10) of the true probability on the $x$ axis, and the error in the bounds on the $y$ axis. Thus, using the solid red line in Figure \ref{fig:modqq}, if the true tail probability of $Q_f$ is $10^{-4}$ ($z=4$), then our new bound for $p=1$ would be approximately of the order $10^{-2}$.

Particularly of note is that while our bounds show an improvement for all functions satisfying the assumptions of Lemma \ref{lemma1}, the improvement is much greater for even functions. This is because our bounds are quadratic, so they must yield the same error bound on both sides of the real line for even functions; however, when bounding an odd function, our bounds will be tight by construction for $x>0$ but may be 
much looser for $x<0$.  As expected, our bounds perform worse for higher powers $p$, which is effectively a result of attempting to control the higher-order behavior of the matrix given traces that measure the empirical mean and variance of the matrix elements.


\section{Conclusions}
We have placed tighter bounds than were previously available on the tails of $Q_f$. Although our bounds are not available in an explicit form, since we optimise over two parameters that previous results set arbitrarily, our bounds are at least as good, which is seen in practice. We further observe that they tend to be significantly tighter and improve relative to the old bounds as we go further out into the tails.

Although our results do give a significantly tighter bound on the tails of $Q_f$, they only work for a specific class of $f$ satisfying the conditions of Lemma 1, which notably excludes functions such as $\exp(x)$. Future developments could improve on this; one possible way would be to introduce an intercept into our quadratic bounds for $\L_1$ and $\L_2$, which would maintain the ease of computability while extending it to a wider range of $f$. A further source of improvement may be achieved by modifying Lemma \ref{lemma1} to account for the asymmetry on $x \leq 0$ vs. $x \geq 0$.  Treating each side of the real line separately could enable one to use both the smallest and largest eigenvalue, rather than just $\underset{i}{\max}\left| \lambda_i\right|$.

Though outside the scope of this paper, it would be possible to achieve similar bounds for sub-Gaussian random variables. This would provide tighter results than currently exist in those cases if the Hanson--Wright inequality argument \cite{rudelson2013hanson} were reworked in terms of explicit constants.

\section{Proofs of Main Results}

\noindent \textbf{Proof of Lemma \ref{lemma1}}

\noindent In the special case $t=0$ we simply have that $\L_1\left(0\right)$, $\L_2(0)$, $\alpha_f(L,0)$, $\beta_f(L,0)$, and $\gamma_f(0)$ are all $0$, so the Lemma clearly holds. We assume now $t \neq 0$.

Since $g_t^1(0)=g_t^2(0)=0$, the choice of $\gamma_f(t)$ is fixed by the need to make 0 a critical point for both of these functions. It remains only to consider the choice of $a$.

Consider $\L$ being either $\L_1$ or $\L_2$.
Write $g(x,a)=\L\left(tf(x)\right) - (ax^2 + bx)$, where $b=\gamma_f(t)$. Since $b$ is fixed, the quadratic functions are
strictly increasing in $a$ at every point. For $x\in (0,L]$ define 
$$
a_x := \frac{\L\left(tf(x)\right)}{x^2} - \frac{\L(tf)'(0)}{x}.
$$ 
Then $a_x$ is the minimum
$a$ such that $g(x,a) \le 0$, and the optimum $a$ that we are looking for is $\sup_{x\in (0,L]} a_x$. We have
\begin{align*}
	\frac{\dif a_x}{\dif x}  & = \frac{\L' \left(h(x)\right)}{x^2} -  \frac{2\L\left(tf(x)\right)}{x^3} + \frac{2 \L(tf)'(0) }{x^2} \\
	  &= 2x^{-3} \left( x\left(\frac{\partial_x \L\left( t f(x)\right)}{2} + t f'(0) \right)  - \L\left(t f(x)\right) \right) \\
	  &\ge 0
\end{align*}
by assumption \ref{E:lemcond1}. Thus $a_x$ is non-decreasing in $x$, and so has its maximum at $L$.
This shows that taking $a=a_L$ makes $g(x,a)\le 0$ for any $x\in (0,L]$, and it is the smallest
such $a$. Note that $a_L= \alpha_f(L)$ when $\L=\L_1$, and $a_L=\beta_f(L)$ when $\L=\L_2$.

Assumption \ref{E:lemcond2} tells us that for $x\in (0,L]$ we have
$\frac{\L\left(tf(x)\right) - \L\left(tf(-x)\right)}{2x} \geq tf'(0)$. This implies that $g(-x,a_L) \le g(x,a_L)\le 0$, so the
same choice of $a=a_L$ provides a bound --- that is, $g(x,a_L) \le 0$ --- over the whole interval $[x,L]$.\qed

\noindent \textbf{Proof of Theorem \ref{Thm1}}

\noindent We credit \cite{Pollard:website} for the proof technique used below.

\noindent Using Lemma 3.1.3 in \cite[p.75]{Christ:2017aa}, for $t < 1/2 d$
\begin{align*}
\mathbb{E}\left[e^{t Q_f }\right] &= \prod\limits_{i=1}^n \left(1-2 t f(\eta_i)  \right)^{-1/2}  \exp\left(\delta_i^2 t f(\eta_i) /(1-2 t f(\eta_i))  \right) \\
&= \exp\left(\sum\limits_{i=1}^n  \delta_i^2 t f(\eta_i)/(1-2 t  f(\eta_i))  - \log\left( 1- 2 t f(\eta_i)\right)/2 \right).
\end{align*}
By Lemma \ref{lemma1} we know, setting $L=\underset{i}{\max}\left| \eta_i\right|$, that for $x\in \left[-L,L\right]$,
\begin{align*}
\L_1\left(tf(x)\right) &\leq \alpha_f(L,t)x^2 + tf'(0)x, \\
\L_2\left(tf(x)\right) & \leq \beta_f(L,t)x^2 + tf'(0)x.
\end{align*}
We claim that this is the optimal choice of $L$. Smaller $L$ will void the inequalities for some $\eta_i$ 
and so cannot be considered. On the other hand, we know that both $\alpha_f(L,t)$ and $\beta_f(L,t)$ are increasing in $L$ so any larger $L$ would simultaneously weaken the quadratic bound 
and shrink the range of values $t$ to which it can be applied, since $1/2f(L)$ is decreasing in $L$.

Therefore,
\begin{align*}
\mathbb{E}\left[e^{tQ_f}\right] &\leq \exp\left(\sum\limits_{i=1}^n  \delta_i^2 \left( \beta_f(L,t)  \eta_i^2 + c \eta_i t \right) +  \alpha_f(L,t) \eta_i^2 + c \eta_i t \right) \\
&\leq \exp\left( \beta_f(L,t) \sum\limits_{i=1}^n  \eta_i^2 \delta_i^2 + c t \sum\limits_{i=1}^n  \eta_i \delta_i^2 + \alpha_f(L,t) \sum\limits_{i=1}^n \eta_i^2 + ct \sum\limits_{i=1}^n  \eta_i \right).
\end{align*}
Applying the definitions of $\xi$ and $\nu_f(t)$ we have
\[
\mathbb{E}\left[\e^{t \left( Q_f - \xi \right) }\right] \leq \e^{\nu_f(t)/2}  .
\]
By Markov's Inequality, for any $q \in \mathbb{R}$,
\begin{multline*}
\mathbb{P}\left( Q_f > q \right) =
\mathbb{P}\left( Q_f - \xi > q -\xi \right)  = \mathbb{P}\left( \e^{Q_f - \xi} > \e^{q -\xi} \right) \\
\leq  \e^{-(q -\xi ) t + \nu_f(t)/2 } \ \text{ for all } t \in \left( 0, 1/2d \right).
\end{multline*}

For $q \leq \xi$, since $\nu_f(t)$ is positive we have the trivial bound $\mathbb{P}\left( Q_f > q \right) \leq 1$.

The bound for $\mathbb{P}\left( Q_f < q \right)$ is derived identically.

\qed

\section*{Acknowledgements}

The first two authors would like to acknowledge the support of the Engineering and Physical Sciences Research Council [grant number EP/M507854/1].  The last author would like to acknowledge the support of the Summer Opportunities Abroad Program (SOAP) - WUSM Global Health \& Medicine and the WUSM Dean's Fellowship, both from the Washington University School of Medicine in St. Louis.

\bibliography{references2}

\end{document}